\documentclass[11pt]{article}
\usepackage{amsmath,amssymb,amsbsy,amsfonts,amsthm,latexsym,
            amsopn,amstext,amsxtra,euscript,amscd,amsthm,url}

\newtheorem{lem}{Lemma}

\newtheorem{thm}{Theorem}

\def\\{\cr}
\def\({\left(}
\def\){\right)}
\def\[{\left[}
\def\]{\right]}
\def\<{\langle}
\def\>{\rangle}

\def\N{{\mathbb N}}
\def\Z{{\mathbb Z}}
\def\K{{\mathbb K}}
\def\Q{{\mathbb Q}}

\newcommand{\bx}{\mathbf{x}}
\newcommand{\dsum}{\displaystyle\sum}

\begin{document}

\title{Some Diophantine equations from finite group theory: $ \Phi_{m}(x) = 2 p^n - 1 $}
\author{Florian Luca (Morelia), Pieter Moree (Bonn) \\ and Benne de Weger (Eindhoven)}
\date{\today}

\pagenumbering{arabic}

\maketitle

\begin{abstract}
We show that the equation in the title (with $ \Phi_n $ the $ n $th cyclotomic polynomial)
has no integer solution with $ n \ge 1 $ in the cases $ (m,p) = (15,41), (15,5581), (10,271) $.
These equations arise in a recent group theoretical investigation by Z.\ Akhlaghi, M.\ Khatami and
B.\ Khosravi.
\end{abstract}
\section{Introduction}
In the recent work \cite{AKK} by Zeinab Akhlaghi, Maryam Khatami and Behrooz Khosravi, some
Diophantine equations come up in a group theoretical context.
In particular, Zeinab Akhlaghi posed the following problems to us.
\begin{itemize}
\item Which primes $ P $ of the form $ P = 2 \cdot 41^{2a} - 1 $ can also be
      written as $ P = \Phi_{15}(q) $, with $ q $ a prime power?
\item Which primes $ P $ of the form $ P = 2 \cdot 5581^{2a} - 1 $ can also be 
      written as $ P= \Phi_{15}(\pm q) $, with $ q $ a prime power?
\item Which primes $P$ of the form $ P=2 \cdot 271^{2a} - 1 $ can also be 
      written as $ P= \Phi_{10}(q^2) $, with $ q $ a prime power?
\end{itemize}
Here $ \Phi_m $ is the $ m $th cyclotomic polynomial. In particular,
\begin{eqnarray*}
   \Phi_{15}(x) & = & x^8 - x^7 + x^5 - x^4 + x^3 - x + 1 , \\
   \Phi_{10}(x) & = & x^4 - x^3 + x^2 - x + 1 .
\end{eqnarray*}
Note that $\Phi_{10}(q^2)=\Phi_{20}(q)$.

Typical for Diophantine equations arising in group theory is the occurrence of primes, and so
the above present some `typical' examples of equations so arising.

Given a group $G$, let $\pi(G)$ denote the set of primes $q$ such that $G$ contains an element
of order $q$. Then the prime graph $\Gamma(G)$ of $G$ is defined as the graph $G$ with vertex
set $\pi(G)$ in which two distinct primes $q,q'\in \pi(G)$ are adjacent if $G$ contains an
element of order $qq'$. Akhlaghi et al.~\cite{AKK} show, using Theorem \ref{thm:main} below and various
already known Diophantine results, that in case $p$ is an odd prime and $k>1$ is odd, then
PGL$(2,p^k)$ is uniqely characterized by its prime graph, i.e. there is no other group having
the same prime graph.

In this paper we will prove the following result, implying the answer "None" to the above three problems,
\begin{thm} \label{thm:main}
Let $ (m,p) = (15,41), (15,5581), $ or $ (10,271) $. Then the Diophantine equation
\begin{equation} \label{eq:main}
   \Phi_m(x) + 1 = 2 p^n
\end{equation}
has no integer solutions $ (n,x) $ with $ n \ge 1 $. 
\end{thm}
In the literature, by different methods, some equations of a similar nature have been 
studied, e.g. the equations  $\Phi_m(x)=p^n$ and $\Phi_m(x)=p^n+1$, with
$m$ a prime. The first equation is a special case of the Nagell-Ljunggren equation
$(x^m-1)/(x-1)=p^n$ and is studied in many papers (for a survey see \cite{BM}). For
a non-existence result of solutions of the second equation see Le \cite{Le}.

General results on solutions of equations of the form $f(x)=by^m$ (see the book
by Shorey and Tijdeman \cite{ST}) imply that for 
an arbitrary, but fixed $m\ge 3$, equation (\ref{eq:main}) has finitely many
solutions $(x,p,n)$ with $\max\{|x|,p,n\}\le C$, with $C$ a computable number. Formulated in
this generality, $C$ which comes from applying the theory of linear forms in logarithms, will
be huge.

In Section \ref{sec:elem}, we give an elementary proof of a lower bound $ n \geq 239 $,
and a related heuristic argument why we do not expect any solutions for these problems.
In Section \ref{sec:bound}, we use algebraic number theory and a deep result
from transcendence theory to deduce an upper bound $ n < 2.163 \cdot 10^{27} $ for
$ n $ satisfying (\ref{eq:main}). Then in Section \ref{sec:lll} the LLL algorithm will be invoked to
efficiently reduce this bound to $ n \leq 59 $. In this way we obtain a rigorous, albeit
computational, proof of Theorem \ref{thm:main}. We note that our method should work in principle
for other equations of the type $ f(x) = a p^n $, when $ f $ is a fixed polynomial with integral
coefficients and at least three distinct roots, $ a \geq 1 $ is a fixed integer, and $ p $ is
a fixed prime not dividing the discriminant of $ f $. The nature of our method is algorithmic
in the sense that for every single choice of parameters the details of the method have to be worked
through separately.

An extended version of this paper including an appendix with numerical material related to
the application of the LLL-method and a subsection on reducing the equation modulo a prime
$q\ne p$ (a subsection not relevant for the proof of our main result), is available as
MPIM preprint 2009-62 \cite{LMW}.

\section{$p$-adic considerations}
\label{sec:elem}
Without loss of generality we may assume that $ |x| \geq 2 $ and $ n \geq 1 $.
We write $ f_m(x) = \Phi_m(x) + 1 $ and $ d = \deg f_m $ for $ m = 10, 15 $.
Elementary calculus shows that for all $ x $
\begin{equation} \label{eq:ineq}
(|x|-1)^d < \Phi_m(x) < f_m(x) < (|x|+1)^d .
\end{equation}
See e.g.~\cite{DT} for some similar estimates.
We start with seeing what information we can derive from studying the $ p $-adic roots of
$ f_m $. If $ (x, n) $ is a solution of (\ref{eq:main}), then there is a root
$$ \bx = \sum_{k=0}^{\infty} a_k p^k \qquad (\textrm{with } a_k \in \{ 0, 1, \ldots, p-1 \}) $$
of $ f_m $ in $ \Q_p $ such that $ x \equiv \bx \pmod{p^n} $. Note that if $ a_0 \neq 0 $ then
the $ p $-adic expansion of $ -\bx $ is
$$ -\bx = (p-a_0) + \sum_{k=1}^{\infty} (p-1-a_k) p^k . $$
Now (\ref{eq:ineq}) with $ f_m(x) = 2 p^n $ implies that
$$ |x| < 2^{1/d} p^{n/d} + 1 < 2 p^{n/d} , $$
and this immediately implies that, in the case $ x > 0 $
$$ a_k = 0 \textrm{ for all } k \in \N \textrm{ with }
   \left\lfloor (n+1)/d \right\rfloor + 1 \leq k \leq n - 1 , $$
and in the case $ x < 0 $
$$ a_k = p-1 \textrm{ for all } k \in \N \textrm{ with }
   \left\lfloor (n+1)/d \right\rfloor + 1 \leq k \leq n - 1 . $$
In other words, the existence of a solution $ n $ of (\ref{eq:main}) implies that of the first
$ n $ $ p $-adic digits of the root $ \bx $, the last consecutive $ \approx n(1-1/d) $ all have
to be equal to $ 0 $ or $ p - 1 $, in respectively the cases $ x > 0 $ and $ x < 0 $.
This seems unlikely to happen, as can easily be verified experimentally for not too large $ n $.
It seems not unreasonable to expect that the $ p $-adic digits of the roots $ \bx $ are uniformly
distributed over $ \{ 0, 1, \ldots, p-1 \} $, and that these distributions per digit are independent.
Then the probability that $n(1-1/d)$ specific consecutive digits are all $ 0 $
respectively $ p - 1 $ is $p^{-n(1-1/d)} $, and the expected number of solutions is at
most $$ \dsum_{n=1}^{\infty} \dfrac{1}{p^{n(1-1/d)}} = \dfrac{1}{p^{(1-1/d)}-1} \ll 1.$$
We conclude that if $p$ is large, then likely there are no solutions. If $p$ is small and there
is no solution with $n$ small, then very likely there are no solutions at all.

A minor variation of the above argument suggests that in case $d\ge 3$ there are only
finitely many solutions $(x,p,n)$ of (\ref{eq:main}) with $n\ge 2$. As we already remarked
in the introduction, this result is known to be true, see \cite{ST}.

Explicit computation of the $ p $-adic root $ \bx $ up to some finite precision
is a quick way to rule out small values of $ n $. We now give details 
for the cases that are of interest to us.
Note that $ p $-adic roots of polynomials are quite easy to compute by Hensel lifting
(i.e.\ the $ p $-adic version of the Newton-Raphson method).

In the case $ p = 41 $ there is one $ 41 $-adic root of $ f_{15} $. Its sequence of
$ 41 $-adic digits is
\begin{eqnarray*}
& ~ & ~ 8, ~18, ~ 3, ~17, ~ 9, ~14, ~12, ~38, ~31, ~35, ~19, ~25, ~19, ~38, ~25, ~24, ~ 1, ~18, \\
& ~ & ~25, ~10, ~14, ~29, ~31, ~18, ~36, ~ 2, ~24, \ldots
\end{eqnarray*}
The smallest $ k $ such that $ a_k = 0 $ or $ 40 $ is $ k = 53 $. Hence a solution of
(\ref{eq:main}) implies $ k \geq 53 $, which in turn implies
$ \left\lfloor (n+1)/8 \right\rfloor + 1 \geq 53 $, so $ n \geq 415 $.

Two remarks are in place. Firstly, we did not even bother to use consecutive zeros, we used only one.
Indeed, $ a_{54} = 15 $, so we could sharpen our result easily.
But we have to stop somewhere, and the result $ n \geq 415 $ is sufficient for the moment.
And secondly, it should be noted that the complexity of this method is exponential, as to compute
the $ n $th $ p $-adic digit we have to compute with numbers of the size $ p^n $. This makes
this method unrealistic for values of $ n $ that become larger than a few thousand.

In the case $ p = 5581 $ there are two $ 5581 $-adic roots of $ f_{15} $. Their sequences of
$ 5581 $-adic digits are
\begin{eqnarray*}
& ~ & ~ 257, ~  64, ~5438, ~1453, ~ 629, ~ 833, ~3090, ~5096, ~4809, ~1493, ~4462, ~1922, \\
& ~ & ~4807, ~ 782, ~3819, ~2190, ~  99, ~2554, ~3603, ~4471, ~1034, ~1407, ~3688, ~\ldots
\end{eqnarray*}
and
\begin{eqnarray*}
& ~ & ~4477, ~3993, ~3590, ~3157, ~3667, ~3404, ~2233, ~3440, ~3784, ~2333, ~ 900, \\
& ~ & ~2522, ~ 184, ~1707, ~5103, ~2005, ~5325, ~1780, ~4765, ~2645, ~3577, ~\ldots
\end{eqnarray*}
In both cases we computed up to $ k = 502 $ and did not encounter a $ 0 $ or a $ 5580 $.
As above it follows that $ n \geq 4015 $.

In the case $ p = 271 $ there is one $ 271 $-adic root of $ f_{10} $. Its sequence of
$ 271 $-adic digits is
\begin{eqnarray*}
& ~& ~241, ~  8, ~147, ~250, ~135, ~263, ~  1, ~126, ~ 89, ~262, ~149, ~ 20, ~147, ~ 78, \\
& ~& ~220, ~219, ~176, ~148, ~206, ~255, ~ 38, ~115, ~186, ~178, ~235, ~\ldots
\end{eqnarray*}
The smallest $ k $ such that $ a_k = 0 $ or $ a_k = 270 $ is $ k = 61 $. Hence a solution of
(\ref{eq:main}) implies $ k \geq 61 $, which in turn implies
$ \left\lfloor (n+1)/4 \right\rfloor + 1 \geq 61 $, so $ n \geq 239 $.

Using the above results we infer that on heuristic grounds with probability at most
$10^{-1000}$ equation (\ref{eq:main}) has a non-trivial solution. Since in mathematics one
has to prove assertions beyond `unreasonable doubt', we cannont conclude our paper
at this point.

\section{Finding an upper bound} \label{sec:bound}
We start with giving some data on relevant algebraic number fields. Then we derive from
equation (\ref{eq:main}) an $ S $-unit inequality, to which we apply transcendence
theory to find an explicit upper bound for $ n $.

\subsection{Field data}
We have
\begin{eqnarray*}
   f_{15}(x) & = & x^8 - x^7 + x^5 - x^4 + x^3 - x + 2 , \\
   f_{10}(x) & = & x^4 - x^3 + x^2 - x + 2 ,
\end{eqnarray*}
and we write equation (\ref{eq:main}) as
\begin{equation} \label{eq:main2}
   f_m(x) = 2 p^n ,
\end{equation}
where $ (m,p) = (15,41), (15,5581), (10,271) $. For brevity we will refer to these cases as
the cases $ p = 41, 5581, 271 $ or $ m = 15, 10 $. In Section \ref{sec:elem}, we have seen
that $ n \geq 239 $, and we may assume that $ |x| \geq 2 $.

The polynomials $ f_m $ are irreducible and have no real roots. We label the roots as follows:
$$ \begin{array}{c|r@{}l|r@{}l}
   \textrm{root} & & f_{15} & & f_{10} \\ \hline
   \alpha^{(1)} & 1.0757\ldots & + 0.4498\ldots i & 0.9734\ldots & + 0.7873\ldots i \\
   \alpha^{(2)} & 0.6243\ldots & + 0.8958\ldots i & -0.4734\ldots & + 1.0255\ldots i \\
   \alpha^{(3)} & -0.1701\ldots & + 1.0292\ldots i &  \\
   \alpha^{(4)} & -1.0299\ldots & + 0.2698\ldots i & & \\ \hline
   & \alpha^{(j)} = \overline{\alpha}^{(j-4)} & \textrm{ for } j=5,6,7,8 &
     \alpha^{(j)} = \overline{\alpha}^{(j-2)} & \textrm{ for } j=3,4 \\ \hline
   \max|\alpha^{(j)}| < & 1.167 & & 1.252 & \\
   \end{array} $$
We write $ \K_m $ for the field $ \Q(\alpha) $ where $ \alpha $ is a root of $ f_m(x) = 0 $,
so that $ d = \deg f_m = [\K_m:\Q] $, i.e.\ $ d = 8 $ for $ m = 15 $ and $ d = 4 $ for $ m = 10 $.

We need a lot of data on these fields. We used Pari \cite{pari} to obtain the data given below.

The discriminants of $ \K_{15}, \K_{10} $ are respectively $ 682862912 = 2^6 \cdot 83 \cdot 128551 $ and
$ 1396 = 2^2 \cdot 349 $. In both cases $ \alpha $ generates a power integral basis. Fundamental units are:
\begin{eqnarray*}
\textrm{for } m = 15 : & & \beta_1 = \alpha^7 + \alpha^4 + \alpha^2 + \alpha - 1, \\
                       & & \beta_2 = \alpha^6 - \alpha^5 + \alpha^4 + \alpha - 1, \\
                       & & \beta_3 = \alpha^2 - \alpha + 1, \\
\textrm{for } m = 10 : & & \beta_1 = \alpha^3 - \alpha^2 + 1.
\end{eqnarray*}
The regulators are $ 4.2219\ldots, 1.1840\ldots $ respectively. The class groups of both
fields is trivial.

The prime decomposition of $ 2 $ is
\begin{eqnarray*}
\textrm{for } m = 15 : & & 2 = \alpha (\alpha+1)^4 (\alpha^3-\alpha^2+1) \beta_1^{-2} \beta_2, \\
\textrm{for } m = 10 : & & 2 = - \alpha (\alpha-1)^3 \beta_1^{-1}.
\end{eqnarray*}
Thus the prime ideals of norm $ 2 $ are $ (\alpha), (\alpha+1) $ when $ m = 15 $, and $ (\alpha), (\alpha-1) $
when $ m = 10 $.

The prime decomposition of $ p $ in the field $ \K_m $ is as follows: \par\noindent
for $ p = 41 $: $ 41 = \gamma_1 \gamma_2 $, where
$$
\begin{array}{ll}
   \gamma_1 = - \alpha^7 + \alpha^6 + \alpha^5 - 2 \alpha^4 + \alpha^3 + \alpha^2 - \alpha + 1, & N(\gamma_1) = 41, \\
   \gamma_2 = - 4 \alpha^7 + 13 \alpha^6 - 19 \alpha^5 + 8 \alpha^4 - 14 \alpha^3 + 7 \alpha^2 + 15 \alpha + 1, &
   N(\gamma_2) = 41^7,
\end{array} $$
for $ p = 5581 $: $ 5581 = \gamma_1 \gamma_2 \gamma_3 \gamma_4 $, where
$$
\begin{array}{ll}
   \gamma_1 = \alpha^6 - \alpha^5 - 2 \alpha + 1, & N(\gamma_1) = 5581, \\
   \gamma_2 = 2 \alpha^5 + \alpha^2 + \alpha + 1, & N(\gamma_2) = 5581, \\
   \gamma_3 = - 3 \alpha^7 - \alpha^6 + 7 \alpha^5 - 4 \alpha^4 - 5 \alpha^3 + 7 \alpha^2 + \alpha + 1, &
   N(\gamma_3) = 5581^2, \\
   \gamma_4 = 85 \alpha^7 - 41 \alpha^6 - 112 \alpha^5 + 55 \alpha^4 - 21 \alpha^3 + & \\
   \hfill + 134 \alpha^2 + 92 \alpha - 135, & N(\gamma_4) = 5581^4,
\end{array} $$
for $ p = 271$: $ 271 = \gamma_1 \gamma_2 $, where
$$
\begin{array}{ll}
   \gamma_1 = - 2 \alpha^3 + 4 \alpha^2 - 4 \alpha + 3, & N(\gamma_1) = 271, \\
   \gamma_2 = - 18 \alpha^3 + 16 \alpha^2 + 44 \alpha + 53, & N(\gamma_2) = 271^3 .
\end{array} $$

\subsection{Deriving an $ S $-unit inequality}
If $x$ is an integer satisfying (\ref{eq:main2}), then it follows that in $ {\mathcal O}_{\K} $ we have
$$ (x-\alpha) z = 2 p^n $$
for a $ z \in {\mathcal O}_{\K} $.
Thus, we can write (taking $ \gamma_3 = \gamma_4 = 1 $ in the cases $ p = 41, 271 $)
$$ x - \alpha = \delta \gamma_1^{n_1} \gamma_2^{n_2} \gamma_3^{n_3} \gamma_4^{n_4} \beta ,
   \qquad z = (2/\delta) \gamma_1^{n-n_1} \gamma_2^{n-n_2} \gamma_3^{n-n_3} \gamma_4^{n-n_4} \beta^{-1} , $$
where $\delta\mid 2$ and $\beta$ is a unit. Taking norms we find
$$ 2 p^n = N(x-\alpha) = N(\delta) p^{c_1n_1+c_2n_2+c_3n_3+c_4n_4} ,$$
where
$ (c_1,c_2,c_3,c_4) = (1,7,0,0), (1,1,2,4), (1,3,0,0) $ for respectively $ p = 41, 5581, 271 $.
It follows that $ N(\delta) = 2 $, and $ n = c_1n_1+c_2n_2+c_3n_3+c_4n_4 $.

First observe that $ 0 < n_i < n $ is impossible. Indeed, for if not, then there exists
$ k \in \{1,2,3,4\} $ such that $ \gamma_k \neq 1 $ divides both $ x - \alpha $ and $ z $. Observe that
if $ \alpha = \alpha^{(i)} $, then $ z = \prod_{j\neq i} (x - \alpha^{(j)}) $. Thus, if $ {\mathfrak p} $
is some prime ideal of $ {\mathcal O}_{\overline{\K}} $ dividing $ \gamma_k $, then $ {\mathfrak p} $
divides both $ x - \alpha^{(i)} $ and $ x - \alpha^{(j)} $ for some $ j \neq i $. In particular,
$ {\mathfrak p} $ divides $ \alpha^{(i)}-\alpha^{(j)} $, and thus also $ \Delta(f_m) $.
Since this last number is an integer and $ {\mathfrak p} $ has norm a power of $ p $
in $ {\overline{\K}_m} $, it would follow that $ p $ divides $ \Delta(f_m) $, which is not the case.
Thus, the only possibilities are $ n_i \in \{ 0, n \} $ for all $ i $. The equation
$ n = c_1 n_1 + c_2 n_2 + c_3 n_3 + c_4 n_4 $ now has only the solutions
$ (n_1,n_2) = (n,0) $ in the cases $ p = 41, 271 $, and $ (n_1,n_2,n_3,n_4) = (n,0,0,0), (0,n,0,0) $
in the case $ p = 5581 $.

We get the following equations:
\begin{equation} \label{eq:1}
\begin{array}{llll}
p = 41:   & x - \alpha = \pm \delta \gamma^n \beta_1^{m_1}\beta_2^{m_2}\beta_3^{m_3}, &
            \delta = \alpha, \alpha+1, & \gamma = \gamma_1, \\
p = 5581: & x - \alpha = \pm \delta \gamma^n \beta_1^{m_1}\beta_2^{m_2}\beta_3^{m_3}, &
            \delta = \alpha, \alpha+1, & \gamma = \gamma_1, \gamma_2 , \\
p = 271:  & x - \alpha = \pm \delta \gamma^n \beta_1^{m_1}, &
            \delta = \alpha, \alpha-1, & \gamma = \gamma_1. \\
\end{array}
\end{equation}
Now we could proceed by conjugating equation \eqref{eq:1} and eliminating $ x $ to get a
unit equation. However, this resulting unit equation will live in the field
$ \Q[\alpha,{\overline{\alpha}}] $, which is of degree $ d(d-1) $, because the Galois group of
$ f_m(x) $ over $ \Q $ is $ S_d $. Since estimates for linear
forms in logarithms are quite sensitive to the degree, we will continue to work in
$ \K_m $. We proceed as follows. For convenience in the cases $ p = 41, 271 $ we put
$ \beta_2 = \beta_3 = 1 $ and $ m_2 = m_3 = 0 $. We have from (\ref{eq:1}) that
\begin{equation} \label{eq:2}
z = \frac{2p^n}{x-\alpha} = \pm \left(\frac{2}{\delta}\right) \left(\frac{p}{\gamma}\right)^n
    \beta_1^{-m_1} \beta_2^{-m_2} \beta_3^{-m_3} .
\end{equation}
Putting $ y = x - \alpha $, Taylor's formula yields
$ z = \dsum_{i=1}^d \frac{f_m^{(i)}(\alpha)}{i!} y^{i-1} $, hence
$$ \left| z - y^{d-1} \right| = \left| \sum_{i=1}^{d-1} \frac{f_m^{(i)}(\alpha)}{i!} y^{i-1} \right| . $$
Let us now make some estimates.
Observe that the lower bound $ n \geq 239 $ from Section \ref{sec:elem} is amply sufficient to
guarantee $ p^n > (2 \cdot 10^{10})^d $. Then (\ref{eq:ineq}) implies
\begin{equation} \label{eq:c1}
|y| = |x-\alpha| \ge |x| - |\alpha| > f_m(x)^{1/d} - 1 - |\alpha| > 2^{1/d} p^{n/d} - 2.252 >
   C_1 p^{n/d} ,
\end{equation}
where $ C_1 = 1.090 $ when $ m = 15 $ and $ C_1 = 1.189 $ when $ m = 10 $.
Hence, $ |y| > 2 \cdot 10^{10} $.
We now compute upper bounds for $ \dfrac{|f_m^{(i)}(\alpha)|}{i!} $, getting
$$
\begin{array}{c|ccccccc}
i & 1 & 2 & 3 & 4 & 5 & 6 & 7 \\ \hline
|f_{15}^{(i)}(\alpha)|/i! < & 16.40 & 56.37 & 109.6 & 126.7 & 90.07 & 39.00 & 9.489 \\
|f_{10}^{(i)}(\alpha)|/i! < & 6.977 & 9.261 & 5.021 & & & & \\
\end{array}
$$
so that
$$ \left| z - y^{d-1} \right| < |y|^{d-2} \sum_{i=1}^{d-1} \frac{f^{(i)}}{i!} \frac{1}{|y|^{d-1-i}} <
   C_2 |y|^{d-2} ,$$
where $ C_2 = 9.490 $ for $ m = 15 $ and $ C_2 = 5.022 $ for $ m = 10 $, because
$ |y| > 2 \cdot 10^{10} $. Thus,
\begin{equation} \label{eq:3}
\left| 1 - \frac{z}{y^{d-1}} \right| < \frac{C_2}{|y|} < {C_3\over p^{n\over d}} ,
\end{equation}
where $ C_3 > \dfrac{C_2}{C_1} $, so $ C_3 = 8.706 $ for $ m = 15 $ and $ C_3 = 4.223 $ for $ m = 10 $.
Using equations \eqref{eq:1}, \eqref{eq:2} and \eqref{eq:3}, we get the $ S $-unit inequality we want:
\begin{equation}
\label{eq:4}
\left|1-\left(\frac{2}{\delta^d}\right)\left(\frac{p}{\gamma^d}\right)^n\beta_1^{-8m_1}\beta_2^{-8m_2}\beta_3^{-8m_3}\right|<\frac{C_3}{p^{n/d}}.
\end{equation}

\subsection{Applying transcendence theory}
We shall apply a linear form in logarithms to bound the expression on the left of inequality \eqref{eq:4}
from below. We first check that it is not zero. If it were, then since it comes from rewriting the left
hand side of inequality \eqref{eq:3}, we would get that $ z = y^{d-1} $. Since $ y z = 2 p^n $, we get that
$ y^d = 2 p^n $, which violates the prime decomposition of $ 2 $ in $ \K_m $.

Next, we need to bound $ m_1, m_2 $ and $ m_3 $ in terms of $ n $.
Since $ p^{n/d} > 2 \cdot 10^{10} $, it follows from (\ref{eq:ineq}) that
\begin{equation} \label{eq:c4}
|y| = |x-\alpha| \le |x| + |\alpha| < f_m(x)^{1/d} + 1 + |\alpha| < 2^{1/d} p^{n/d} + 2.252 <
   C_4 p^{n/d} ,
\end{equation}
where $ C_4 = 1.091 $ for $ m = 15 $ and $ C_4 = 1.190 $ for $ m = 10 $.
Now we take absolute values of the conjugates of equation \eqref{eq:1}, and rewrite them as
\begin{equation} \label{eq:i}
\frac{\left|x-\alpha^{(i)}\right|}{\left|\delta^{(i)}\right|\left|\gamma_1^{(i)}\right|^n} =
\left|\beta_1^{(i)}\right|^{m_1} \left|\beta_2^{(i)}\right|^{m_2} \left|\beta_3^{(i)}\right|^{m_3} .
\end{equation}
We computed:
$$
\begin{array}{ccc}
\textrm{for } p = 41:   & 0.2714 < |\delta^{(i)}| < 2.124 , & 0.5676 < |\gamma^{(i)}| < 5.349 , \\
\textrm{for } p = 5581: & 0.2714 < |\delta^{(i)}| < 2.124 , & 1.522  < |\gamma^{(i)}| < 5.531 , \\
\textrm{for } p = 271:  & 0.7877 < |\delta^{(i)}| < 1.796 ,&  2.253  < |\gamma^{(i)}| < 7.307 , \\
\end{array}
$$
and thus
\begin{eqnarray*}
& & \max\(\log\frac{p^{1/d}}{\min|\gamma^{(i)}|},\log\frac{\max|\gamma^{(i)}|}{p^{1/d}}\) < C_5, \\
& & \max\(\log\frac{C_4}{\min|\delta^{(i)}|},\log\frac{\max|\delta^{(i)}|}{C_1}\) < C_6 ,
\end{eqnarray*}
where for $ p = 41 $ we have $ C_5 = 1.213, C_6 = 1.392 $, for $ p = 5581 $ we have $ C_5 = 0.6584, C_6 = 1.392 $, and for
$ p = 271 $ we have $ C_5 = 0.5884, C_6 = 0.4126 $.
It follows from (\ref{eq:c1}) and (\ref{eq:c4}) that
$$ \left|\log\left(\frac{\left|x-\alpha^{(i)}\right|}{
   \left|\delta^{(i)}\right|\left|\gamma^{(i)}\right|^n}\right)\right| < C_5 n + C_6. $$
Writing $ u_i $ for the logarithm of the left hand side of equations \eqref{eq:i}, we get that
\begin{equation} \label{eq:sys}
u_i = m_1 \log\left|\beta_1^{(i)}\right| + m_2 \log\left|\beta_2^{(i)}\right| +
      m_3 \log\left|\beta_3^{(i)}\right| \textrm{ for three conjugates } i,
\end{equation}
and hence $ |u_i| < C_5 n + C_6 $ for all $ i $. If $ m = 10 $ this simply states
$ \log\left|\beta_1^{(i)}\right| > 1.184 $ (this is the regulator of $ K_{10} $), as then
$ \beta_2 = \beta_3 = 1 $,
and thus $ |m| < (C_5 n + C_6)/1.184 $. If $ m = 15 $, solving the system (\ref{eq:sys})
with Cramer's rule, we get that
$$ \max\{|m_1|,|m_2|,|m_3|\}<\frac{3(C_5n+C_6)R_2}{R_{\beta}}, $$
where $ R_2 $ is the maximal absolute value of all the $ 2\times 2 $ minors of the coefficient matrix
appearing in formula \eqref{eq:sys} whose determinant is $ R_{\beta} $. The minor largest in absolute value is the $(2,1)$ minor obtained by eliminating the second row and first column, and its value is
$ R_2 < 2.746 $. Putting all this together gives
$$ \max\{|m_1|,|m_2|,|m_3|\} < C_7 n + C_8, $$
where $ C_7 = 2.369, C_8 = 2.718 $ when $ p = 41 $, $ C_7 = 1.286, C_8 = 2.718 $ when $ p = 5581 $, and
$ C_7 = 0.4970, C_8 = 0.3485 $ when $ p = 271 $.

The next step is to prepare for the application of a deep result from transcendence theory.
We return to inequality \eqref{eq:4} and rewrite it as
\begin{equation} \label{eq:imp}
\left| 1 - \prod_{i=1}^r \eta_i^{b_i} \right| < \frac{C_3}{p^{n/d}},
\end{equation}
where $ r = 5 $ when $ m = 15 $ and $ r = 3 $ if $ m = 10 $, and
$$ \eta_1 = \frac{2}{\delta^d}, \quad \eta_2 = \frac{p}{\gamma^d}, \quad \eta_3 = \beta_1, \quad
   \eta_4 = \beta_2, \quad \eta_5 = \beta_3, $$
and $ b_1 = 1, b_2 = n, b_3 = -d m_1, b_4 = -d m_2, b_5 = -d m_3 $ are integers satisfying
$$ B = \max|b_i| < d (C_7 n + C_8). $$
Recall that for an algebraic number $\eta$ having 
$$ a_0 \prod_{i=1}^d (X-\eta^{(i)})$$
as minimal polynomial over the integers,
the logarithmic height is defined as
$$ h(\eta) = \frac{1}{d} \left( \log|a_0| + \sum_{i=1}^d \log \max \left\{ \left|\eta^{(i)}\right|,
   1 \right\} \right). $$
With this notation, Matveev \cite{Mat} proved the following deep theorem.

\begin{thm} \label{thm:Mat}
Let $ \K $ be a field of degree $ D $, $ \eta_1, \ldots, \eta_k $ be nonzero elements of $ \K $, and
$ b_1, \ldots, b_k $ integers. Put
$$ B = \max\{|b_1|, \ldots, |b_k|\} $$
and
$$ \Lambda = 1 - \prod_{i=1}^k \eta_i^{b_i}. $$
Let $ A_1, \ldots, A_k $ be real numbers such that
$$ A_j \ge \max\{ D h(\eta_j), |\log \eta_j|, 0.16\}, \qquad j = 1, \ldots, k. $$
Then, assuming that $ \Lambda \ne 0 $, we have
$$ \log |\Lambda| > -3 \cdot 30^{k+4} (k+1)^{5.5} D^2 (1+\log D) (1+\log(kB)) \prod_{i=1}^k A_i. $$
\end{thm}

We apply Matveev's result to get a lower bound on the expression appearing in the left hand side of \eqref{eq:imp}  with $k=r+1$. We take the field to be our $\K_m$, so $D=d$. We also take
$ \eta_i, b_i $ as in \eqref{eq:imp}.

We computed as leading coefficients $ a_0 $ of minimal polynomials:
$$
\begin{array}{r|l|lll}
 m & \delta   & \eta_1    & \eta_2    & \eta_3, \eta_4, \eta_5 \\ \hline
15 & \alpha   & a_0 = 2^7 & a_0 = p^7 & a_0 = 1                \\
   & \alpha+1 & a_0 = 2^4 &           &                        \\ \hline
10 & \alpha   & a_0 = 2^3 & a_0 = p^3 & a_0 = 1                \\
   & \alpha-1 & a_0 = 2   &           &                        \\
\end{array}
$$
and for the $ A_j $ we found
$$
\begin{array}{r|cccccc}
   p & A_1 < & A_2 < & A_3 < & A_4 < & A_5 < \\ \hline
  41 & 25.02 & 47.80 & 4.371 & 4.247 & 2.976 \\
5581 & 25.02 & 74.22 & 4.371 & 4.247 & 2.976 \\
 271 & 3.988 & 21.52 & 2.634 &  &  \\
\end{array}
$$
Thus, by Matveev's bound we have that
$$ |\log \Lambda| > -C_9 (1 + \log(rB)), $$
where $ C_9 > 3 \cdot 30^{r+4} (r+1)^{5.5} d^2 (1+\log d) A_1 A_2 \ldots A_r $ satisfies
\begin{eqnarray*}
\textrm{for } p = 41:   & & C_9 = 1.465 \cdot 10^{25}, \\
\textrm{for } p = 5581: & & C_9 = 2.275 \cdot 10^{25}, \\
\textrm{for } p = 271:  & & C_9 = 1.160 \cdot 10^{18}.
\end{eqnarray*}
Comparing this with the fact that $B\le d(C_7n+C_8)$ and with inequality \eqref{eq:imp}, we get
$$ \frac{\log p}{d} n - \log C_3 < - \log|\Lambda| < C_9 (1 + \log(rd(C_7n+C_8))). $$
Concretely:
\begin{eqnarray*}
\textrm{for } p = 41:   & & 0.4641 n - 2.165 < 1.465 \cdot 10^{25} (1+\log(94.80n+108.8)) \\
                        & & \textrm{implying } n < N = 2.163 \cdot 10^{27}, \\
\textrm{for } p = 5581: & & 1.078 n - 2.165 < 2.275 \cdot 10^{25} (1+\log(51.45n+108.8)) \\
                        & & \textrm{implying } n < N = 1.424 \cdot 10^{27}, \\
\textrm{for } p = 271:  & & 1.400 n - 1.441 < 1.160 \cdot 10^{18} (1+\log(5.964n+4.182)) \\
                        & & \textrm{implying } n < N = 3.970 \cdot 10^{19}.
\end{eqnarray*}

\section{Reducing the upper bound} \label{sec:lll}
So, it remains to solve
$$ \left\{ \begin{array}{rcl}
   \left|1-\eta_1\eta_2^{n} \displaystyle \prod_{i=3}^r \eta_i^{-dm_{i-2}}\right| & < & \dfrac{C_3}{p^{n/d}} , \\
   \max|m_i| & < & d(C_7n+C_8) , \\
   n & < & N.
   \end{array} \right. $$
This is a finite problem, but the upper bound $ N $ is way too large to apply
brute force or the method from Section \ref{sec:elem}. Efficient methods for solving such
problems based on lattice basis reduction using the LLL algorithm exist, see \cite{dW}, and
they work quite well in our case. Here are the details.

We put
$$ \lambda_i^{(j)} = \left\{ \begin{array}{rcl}
   \log\left|\eta_i^{(j)}\right| & {\text{\rm for}} & i = 1, 2 \\
-d \log\left|\eta_i^{(j)}\right| & {\text{\rm for}} & i = 3, \ldots, r.
\end{array} \right., \quad j = 1, \ldots, r-1 . $$
Let
$$ \lambda^{(j)} = \lambda_1^{(j)} + n \lambda_2^{(j)} + m_1 \lambda_3^{(j)} + \ldots +
   m_{r-2} \lambda_r^{(j)} \quad {\text{\rm for}}~j = 1, \ldots, r-1 . $$
By (\ref{eq:imp}), the real linear forms $ \lambda^{(j)} $ satisfy
\begin{equation} \label{eq:lam}
\left|\lambda^{(j)}\right| \le \left|1-e^{\lambda^{(j)}}\right| \leq \left|1-\eta_1^{(j)} \left(\eta_2^{(j)}\right)^{n} \prod_{i=3}^r \left(\eta_i^{(j)}\right)^{-dm_{i-2}} \right| < \frac{C_3}{p^{n/d}} .
\end{equation}
We let $ K $ be some constant slightly larger than $ N^{(r-1)/(r-2)} $, i.e.\ $ N^{4/3} $
when $ m = 15 $ and $ r = 5 $, and $ N^2 $ when $ m = 10 $ and $ r = 3 $. We write
$ \theta_i^{(j)} = \left[ K \lambda_i^{(j)} \right] $ for $ i = 1, \ldots, r $, where $ [\cdot] $
denotes rounding to the nearest integer. We put
$$ (\lambda')^{(j)} = \theta_1^{(j)} + n \theta_2^{(j)} + m_1 \theta_3^{(j)} + \ldots +
   m_{r-2} \theta_r^{(j)} . $$
Then
$$ \left| K \lambda^{(j)} - (\lambda')^{(j)}\right| \le \frac{1}{2} + \frac{n}{2} +
   \frac{r-2}{2} \max|m_i| < C_{10} n + C_{11} , $$
where $ C_{10} = \frac{1}{2} + \frac{r-2}{2} d C_7 $ and
$ C_{11} = \frac{1}{2} + \frac{r-2}{2} d C_8 $.
Then $ n \geq N $ implies
\begin{equation} \label{eq:lamaks}
\left|(\lambda')^{(j)}\right| < K \left| \lambda^{(j)} \right| + C_{10} N + C_{11} .
\end{equation}
We now look at the matrix $\Gamma$ and the vector ${\underline{y}}$ given as
\begin{eqnarray*}
\textrm{for } m = 15: & & \Gamma = \left( \begin{matrix}
\theta_3^{(i)} & \theta_4^{(i)} & \theta_5^{(i)} & \theta_2^{(i)} \\
\theta_3^{(j)} & \theta_4^{(j)} & \theta_5^{(j)} & \theta_2^{(j)} \\
\theta_3^{(k)} & \theta_4^{(k)} & \theta_5^{(k)} & \theta_2^{(k)} \\
0 & 0 & 0 & 1
\end{matrix} \right), \quad {\underline{y}} = \left( \begin{matrix}
-\theta_1^{(i)} \\ -\theta_1^{(j)} \\ -\theta_1^{(k)} \\ 0
\end{matrix} \right), \\
& & \textrm{where } (i,j,k)\in \{(1,2,3),(1,2,4),(1,3,4),(2,3,4)\}, \\
\textrm{for } m = 10: & & \Gamma = \left( \begin{matrix}
\theta_3^{(i)} & \theta_2^{(i)} \\ 0 & 1
\end{matrix} \right), \quad {\underline{y}} = \left( \begin{matrix}
-\theta_1^{(i)} \\ 0
\end{matrix} \right), \\
& & \textrm{where } i \in \{1,2\}.
\end{eqnarray*}
Observe that for $ {\underline{x}} = (m_1,\ldots,m_{r-2},n)^T $
$$ \Gamma{\underline{x}}-{\underline{y}} = \left((\lambda')^{(i)}, (\lambda')^{(j)}, (\lambda')^{(k)},
   n \right)^T \quad \textrm{ resp.\ } \quad \left((\lambda')^{(i)}, n \right)^T . $$
The columns of $\Gamma$ generate a sublattice of $\Z^{r-2}$. Let
$ d(\Gamma,{\underline{y}}) = \displaystyle\min_{{\underline{x}} \in \Z^{r-2}}
\left|\Gamma{\underline{x}}-{\underline{y}}\right| $
be the distance from $ {\underline{y}} $ to the nearest lattice point. From (\ref{eq:lamaks}) we find
\begin{equation} \label{eq:d}
d(\Gamma,{\underline{y}}) \le \left| \Gamma {\underline{x}}-{\underline{y}} \right| <
\sqrt{(r-2) \left( K \max \left|\lambda^{(j)}\right| + C_{10} N + C_{11} \right)^2 + N^2}.
\end{equation}
Put
$$ c = \frac{N^{1/(r-2)}}{K} \left( {\sqrt{\frac{d(\Gamma,{\underline{y}})^2-N^2}{r-2}}}
       - (C_{10}N+C_{11}) \right). $$
If $c$ happens to be a positive real number, then combining (\ref{eq:lam}) and (\ref{eq:d}) we get
for $ \lambda = \lambda^{(j)}$, such that $ |\lambda| = \max \left|\lambda^{(j)}\right| $ satisfies
$$ c N^{-1/(r-2)} < |\lambda| < \frac{C_3}{p^{n/d}}, $$
and hence
$$ n < \frac{d}{\log p} \left(\log C_3 - \log c + \frac{1}{r-2} \log N \right). $$
In particular, if $c$ is reasonable, that is, not too tiny, then the above bound is a reduced upper
bound for $ n $. We can argue that this is reasonable, because if the lattice is generic, that is,
if it satisfies
$$ d(\Gamma, {\underline{y}}) \approx {\text{\rm det}}(\Gamma)^{1/\dim \Gamma} \approx K^{(r-2)/(r-1)}, $$
then with the choice of $ K $ being somewhat larger than $ N^{(r-1)/(r-2)} $, one would expect that
$ d(\Gamma,{\underline{y}}) $ is somewhat larger than $ N $, so that $ c $ just becomes positive:
$$ c \approx \frac{N^{1/(r-2)}}{K} \cdot N \approx 1. $$
Clearly, a lower bound for $ d(\Gamma,{\underline{y}}) $ suffices. To compute such a bound we use
Lemma 3.5 from \cite{dW}, which we now state.

\begin{lem}
If $ {\underline{c}_1}, \ldots, {\underline{c}_{r-1}} $ is an LLL-reduced basis for the lattice
spanned by the columns of the matrix $ \Gamma $, and $ (s_1, \ldots, s_{r-1}) $ are the coordinates of
$ {\underline{y}} \in \Z^{r-1} $ with respect to this basis, then
$$ d(\Gamma,{\underline{y}}) \ge 2^{-(r-2)/2} \|s_{r-1}\| |{\underline{c}_1}|, $$
where $ \|\cdot\| $ denotes the distance to the nearest integer.
\end{lem}

When a new upper $ N_1 $ on $ n $ is found, the procedure can be repeated with $ N_1 $ instead of $ N $.

As for the practical calculations, for $ p = 41 $ and $ p = 5581 $ we use $ K = 10^{39} $,
and for $ p = 271 $ we use $ K = 10^{41} $.
For $ p = 41 $ the conjugates $ (i,j,k) = (1,3,4) $ turned out to give the best results,
and for $ p = 5581 $ we took the conjugates $ (i,j,k) = (2,3,4) $ in the case $ \gamma = \gamma_1 $,
and $ (i,j,k) = (1,3,4) $ in the case $ \gamma = \gamma_2 $. For $ p = 271 $ we took
the conjugate $ i = 2 $.
The values of the entries of $ \Gamma $ and $ \underline{y} $ are given in the appendix of
the extended version of this paper \cite{LMW}.

As a result of our computations we found: \medskip \par\noindent
\underline{for $ p = 41 $}: $ | \underline{c}_1| = 1.148\ldots \cdot 10^{30} $, \par\noindent
\begin{tabular}{llll}
for $ \delta = \alpha $: & $ \| s_4 \| = 0.2505\ldots $, & $ d(\Gamma,\underline{y})
\geq 1.017 \cdot 10^{29} $, & $ c = 0.0650\ldots $, \\
for $ \delta = \alpha + 1$: & $ \| s_4 \| = 0.0809\ldots $, & $ d(\Gamma,\underline{y})
\geq 3.286 \cdot 10^{28} $, & $ c = 0.0125\ldots $. \end{tabular}
We infer $ n \leq N_1 = 59 $. \medskip \par\noindent
\underline{for $ p = 5581, \gamma = \gamma_1 $}: $ | \underline{c}_1| = 1.123\ldots \cdot 10^{30} $,
\par\noindent
\begin{tabular}{llll}
for $ \delta = \alpha $: & $ \| s_4 \| = 0.4489\ldots $, & $ d(\Gamma,\underline{y})
\geq 1.784 \cdot 10^{29} $, & $ c = 0.1119\ldots $, \\
for $ \delta = \alpha + 1$: & $ \| s_4 \| = 0.3512\ldots $, & $ d(\Gamma,\underline{y})
\geq 1.395 \cdot 10^{29} $, & $ c = 0.0867\ldots $. \end{tabular}
We infer $ n \leq N_1 = 23 $. \medskip \par\noindent
\underline{for $ p = 5581, \gamma = \gamma_2 $}: $ | \underline{c}_1| = 6.875\ldots \cdot 10^{29} $,
\par\noindent
\begin{tabular}{llll}
for $ \delta = \alpha $: & $ \| s_4 \| = 0.3849\ldots $, & $ d(\Gamma,\underline{y})
\geq 9.357 \cdot 10^{28} $, & $ c = 0.0568\ldots $, \\
for $ \delta = \alpha + 1$: & $ \| s_4 \| = 0.4225\ldots $, & $ d(\Gamma,\underline{y})
\geq 1.027 \cdot 10^{29} $, & $ c = 0.0628\ldots $. \end{tabular}
We infer $ n \leq N_1 = 23 $. \medskip \par\noindent
\underline{for $ p = 271 $}: $ | \underline{c}_1| = 2.826\ldots \cdot 10^{20} $, \par\noindent
\begin{tabular}{llll}
for $ \delta = \alpha $: & $ \| s_2 \| = 0.2302\ldots $, & $ d(\Gamma,\underline{y})
\geq 4.602 \cdot 10^{19} $, & $ c = 0.0014\ldots $, \\
for $ \delta = \alpha - 1$: & $ \| s_2 \| = 0.2565\ldots $, & $ d(\Gamma,\underline{y})
\geq 5.127 \cdot 10^{19} $, & $ c = 0.0050\ldots $. \end{tabular}
We infer $ n \leq N_1 = 37 $. \medskip

All reduced upper bounds are well below the lower bound $ n \geq 239 $ we had already
found in Section \ref{sec:elem}. Hence, the given equations have no positive integer solutions $ (n,x) $.

We used the built-in LLL implementation of Mathematica 7.0.
The total computation time was about 0.5 second on a standard laptop. \hfil\break

{\tt Acknowledgement}. We like to thank both N.~Bruin and S.~Akhtari
for sketching other approaches to proving Theorem \ref{thm:main}. Both these
approaches seem more involved than ours. On the other hand, we cannot exclude
that less preliminary considerations on their part would lead to a shorter
proof than ours. Also, we thank M.~Bennett for some helpful remarks.

The information on value sets given in Section 2.2 of the extended
version \cite{LMW} was kindly provided to
us by D.~Wan and N.~Alexander.

Work on this paper started during a visit of the first author to the Max-Planck-Institute
of Mathematics in the Spring of 2009.

\vfil\eject

\medskip\noindent Florian Luca \par\noindent
{\footnotesize {Instituto de Matem{\'a}ticas},{ Universidad Nacional Autonoma de M{\'e}xico} \hfil\break
{C.P. 58089, Morelia, Michoac{\'a}n, M{\'e}xico} \hfil\break
e-mail: {\tt fluca@matmor.unam.mx}}

\medskip\noindent Pieter Moree \par\noindent
{\footnotesize Max-Planck-Institut f\"ur Mathematik,
Vivatsgasse 7, D-53111 Bonn, Germany.\hfil\break
e-mail: {\tt moree@mpim-bonn.mpg.de}}

\medskip\noindent Benne de Weger \par\noindent
{\footnotesize
Department of Mathematics and Computer Science,\hfil\break
Eindhoven University of Technology,\hfil\break
PO Box 513, 5600 MB Eindhoven, The Netherlands,\hfil\break
e-mail: {\tt b.m.m.d.weger@tue.nl}}


\begin{thebibliography}{9999}

\bibitem{AKK} Z.~Akhlaghi, M.~Khatami and B.~Khosravi,
Characterization of $PGL(2,p^k)$, by prime graphs, where $k>1$, is odd, submitted for publication.

\bibitem{BM} Y.~Bugeaud and M.~Mignotte, L'\'equation de
Nagell-Ljunggren $\frac{x\sp n-1}{x-1}=y\sp q$, {\it Enseign. Math.} (2)
{\bf 48} (2002), 147--168.

\bibitem{DT} G.~Drauschke and M.~Tasche, Prime factorizations for values of cyclotomic 
polynomials in ${\mathbb Z}[i]$, {\it Arch. Math. (Basel)}  {\bf 49}  (1987),  292--300.

\bibitem{Le} M.H.~Le, A note on the
Diophantine equation $(x\sp m-1)/(x-1)=y\sp n+1$,
{\it Math. Proc. Cambridge Philos. Soc.} {\bf 116} (1994), 385--389.

\bibitem{LMW} F.~Luca, P.~Moree and B.~de Weger, Some Diophantine equations from finite group 
theory: $ \Phi_{m}(x) = 2 p^n - 1$, MPIM-preprint 2009-62, \hfil\break
http://www.mpim-bonn.mpg.de/Research/MPIM+Preprint+Series/


\bibitem{Mat} E.~M.~Matveev, An explicit lower bound for a homogeneous rational linear form in logarithms of algebraic numbers. II. (Russian)
{\it Izv. Ross. Akad. Nauk Ser. Mat.\/} {\bf 64}  (2000),  no. 6, 125--180;  translation in {\it Izv. Math.\/} {\bf  64}  (2000),  no. 6, 1217--1269.

\bibitem{pari}
    PARI/GP, version {\tt 2.3.4}, Bordeaux, 2006, \url{http://pari.math.u-bordeaux.fr/}.
    
\bibitem{ST}   T.N.~Shorey and R.~Tijdeman, {\it Exponential Diophantine equations}, 
Cambridge Tracts in Mathematics {\bf 87}, Cambridge University Press, Cambridge, 1986.

\bibitem{dW} B.~M.~M.~de Weger, {\it Algorithms for Diophantine equations\/}, CWI Tract {\bf 65}, Stichting Mathematisch Centrum,
Centrum voor Wiskunde en Informatica, Amsterdam, 1989.

\end{thebibliography}
\end{document}